\documentclass[10pt]{amsart}
\usepackage[cp1251]{inputenc}
\usepackage[english,russian]{babel}
\usepackage{amsmath}
\usepackage{amssymb}
\usepackage{amsfonts}

\newtheorem{theorem}{Theorem}

\newtheorem{corollary}{Corollary}

\newtheorem{question}{Question}
\newtheorem{problem}{Problem}

\def\Z{\mathbb Z}

\def \pw{{\rm pw}}

\begin{document}
\renewcommand{\refname}{References}
\renewcommand{\proofname}{Proof.}
\thispagestyle{empty}

\title[Some Problems]{Some Problems on Knots, Braids, and Automorphism Groups}

\author{{V.~Bardakov}}%
\address{Valery Georgievich Bardakov
\newline\hphantom{iii} Sobolev Institute of Mathematics, Novosibirsk State University pr. Koptyuga, 4, 630090, Novosibirsk, Russia}
\email{bardakov@math.nsc.ru}%

\author{K. Gongopadhyay}%
\address{Krishnendu Gongopadhyay
\newline\hphantom{iii} Department of Mathematical Sciences, Indian Institute of Science Education and Research (IISER) Mohali, Knowledge City, Sector 81, S.A.S. Nagar, P.O. Manauli 140306, India}
\email{krishnendu@iisermohali.ac.in, krishnendug@gmail.com}%

\author{M. Singh}%
\address{Mahender Singh
\newline\hphantom{iii} Department of Mathematical Sciences, Indian Institute of Science Education and Research (IISER) Mohali, Knowledge City, Sector 81, S.A.S. Nagar, P.O. Manauli 140306, India}
\email{mahender@iisermohali.ac.in}%

\author{A.~Vesnin}%
\address{Andrei Vesnin
\newline\hphantom{iii} Sobolev Institute of Mathematics, Novosibirsk State University pr. Koptyuga, 4, 630090, Novosibirsk, Russia}
\email{vesnin@math.nsc.ru}%

\author{J. Wu}%
\address{Jie Wu
\newline\hphantom{iii} Department of Mathematics, National University of Singapore, 10 lower Kent Ridge Road,  119076, Singapore}
 \email{matwuj@nus.edu.sg}%

\thanks{\sc Bardakov, V.G., Gongopadhyay, K., Singh, M., Vesnin, A., Wu, J.,
Some problems on knots, braids, and automorphism groups}
\thanks{\copyright \ 2015 Bardakov, V.G., Gongopadhyay, K., Singh, M., Vesnin, A., Wu, J.}
\thanks{\rm The work is supported by the Indo-Russian DST-RFBR project grant DST/INT/RFBR/P-137 and RFBR-13-01-92697.}
\maketitle
\begin{quote}
\noindent{\sc Abstract. }
We present and discuss some open problems formulated by participants of the International Workshop ``Knots, Braids, and Auto\-mor\-phism Groups'' held in Novosibirsk, 2014.
Problems are related to palindromic and commutator widths of groups; properties of Brunnian braids and two-colored braids, corresponding to an amalga\-ma\-tion of groups; extreme properties of hyperbolic 3-orbifold groups, relations bet\-ween inner and quasi-inner automorphisms of groups; and Staic's construc\-tion of symmetric cohomology of groups. 

\medskip

 \end{quote}

\date{\today}

\section*{Introduction} \label{sec0}
The aim of this paper is to present some problems formulated by participants of the International Workshop ``Knots, Braids and Automorphism Groups'' held at Sobolev Institute of Mathematics (Novosibirsk, Russia) in July, 21--25, 2014. The list of participants, the program, and abstracts of talks are available on the workshop homepage \cite{page}.  

Various open problems were presented and discussed during the problem session. They were motivated by talks given during the workshop as well as research interests of the participants. We present some of those problems in this paper. For some problems, we give preliminary discussions and useful references. In section~\ref{sec1}, we discuss palindromic and commutator widths of groups, and give estimates of these widths for some classes of groups. In section~\ref{sec2}, we formulate problems on braids with some special properties: Brunnian braids and two-colored braids, corresponding to an amalga\-ma\-tion of groups. In section~\ref{sec3}, we discuss intersecting subgroups of link groups, extreme properties of hyperbolic 3-orbifold groups, linearity of 3-manifold groups, and relations between inner and quasi-inner automorphisms of groups. In section~\ref{sec4}, we present Staic's construction of symmetric cohomology of groups and formulate some problems related to it. 

\section{Palindromic widths of groups} \label{sec1}

Let $G$ be a group with a set of generators $X$. A reduced word in the alphabet $X^{\pm 1}$ is a \emph{palindrome} if it reads the same forwards and backwards. The palindromic length $l_{\mathcal P}(g)$ of an element $g$ in $G$ is the minimum number $k$ such that $g$ can be expressed as a product of $k$ palindromes. The \emph{palindromic width} of $G$ with respect to $X$ is defined to be
$$
{\rm pw}(G, X) = \underset{g \in G}{\sup} \ l_{\mathcal{P}}(g).
$$
When there is no confusion about the underlying generating set $X$, we simply denote the palindromic width with respect to $X$ by ${\rm pw}(G)$. Palindromes in groups have already proved useful in studying various aspects of combinatorial group theory and geometry, for example see  \cite{BST}--\cite{P}. It was proved in \cite{BST} that the palindromic width of  a non-abelian free group is infinite.  This result was generalized in \cite{BT} where the authors proved that almost all free products have infinite palindromic width; the only exception is given by the free product of two cyclic groups of order two, when the palindromic width is two. Piggot \cite{P} studied the relationship between primitive words and palindromes in free groups of rank two. It follows from \cite{BST, P} that up to conjugacy, a primitive word can always be written as either a palindrome or a product of two palindromes and that certain pairs of palindromes will generate the group.

Bardakov and Gongopadhyay \cite{BG} initiated the investigation of palindromic width of a finitely generated group that is free in some variety of groups. They demonstrated finiteness of palindromic widths of free nilpotent and certain solvable groups. The following results were established.

\begin{theorem}\label{bgn} \cite{BG, BG2}
Let ${\rm N}_{n, r} $ be the free $r$-step nilpotent group of rank $n \geq 2$. Then the  following holds:
\begin{enumerate}
\item{ The palindromic width ${\rm pw}({\rm N}_{n, 1})$ of a free abelian group of rank $n$  is equal to $n$. }
\item{For $r \geq 2$, $2(n-1) \leq {\rm pw}({\rm N}_{n, r} ) \leq 3n$. }
\item{$2(n-1)\leq {\rm pw}({\rm N}_{n, 2}) \leq 3 (n-1)$. }
\end{enumerate}
\end{theorem}

\begin{problem}{\rm (V. Bardakov~--- K. Gongopadhyay)} \label{prob1}
\begin{enumerate}
\item{For $n \geq 3$, $r \geq 2$, find ${\rm pw}({\rm N}_{n,r} )$.}
\item{Construct an algorithm that determines $l_{\mathcal P} (g)$ for arbitrary $g \in {\rm N}_{n, r} $.}
\end{enumerate}
\end{problem}

We recall that a group $G$ is said to satisfy  \emph{the maximal condition for normal subgroups} if every normal subgroup of $G$ is the normal closure of a finite subset of~$G$. 

\newpage
\begin{theorem}\label{bgs}  \cite{BG3}
\begin{enumerate}
\item {Let $A$ be a normal abelian subgroup of a finitely generated solvable group $G=\langle X \rangle$ such that $G/A$ satisfies the maximal condition for normal subgroups. Then $\pw(G, X) < \infty$.}
\item{ Let $G$ be a finitely generated free abelian-by-nilpotent-by-nilpotent group. Then the palindromic width of $G$ is finite.}
\end{enumerate}
\end{theorem}

As a corollary to (2) of the above theorem the following result holds.

\begin{corollary}\label{sol3} \cite{BG3}
Every finitely generated $3$-step solvable group has finite palindro\-mic width with respect to any finite generating set.
\end{corollary}

Riley and Sale \cite{RS} studied palindromic width of solvable and metabelian groups using different techniques. They proved the special case of Theorem~\ref{bgs}(1) when $A$ is the trivial subgroup: the palindromic width of a finitely generated solvable group that satisfies the maximal condition for normal subgroups has finite palindromic width. All mention results concerning solvable groups established finiteness of pa\-lin\-dro\-mic widths in these groups. However, the techniques that have been used to prove these results do not provide a bound for the widths. This poses the following problem.

\begin{problem}[V. Bardakov~--- K. Gongopadhyay]
Let $G$ be a group in  Theorem \ref{bgs}. Find the lower and upper bounds on the palindromic width $\pw(G)$.
\end{problem}

It would be interesting to understand palindromic width of wreath products. It was proved in \cite{BG3} and \cite{RS} independently  that ${\rm pw}(\Z \wr \Z)=3$. The proof from \cite{RS} is based on the estimate of $\pw(G \wr \Z^r)$, where $G$ is a finitely generated group. The proof from \cite{BG3} relies on the fact that any element in the commutator subgroup of $\Z \wr \Z$ is a commutator. Palindromic width in wreath products has also been investigated by Fink \cite{fink} who has also obtained an estimate of $\pw(G\wr \Z^r)$. All these results are related to the following general problem.

\begin{problem}[V. Bardakov~--- K. Gongopadhyay]
Let $G = A \wr B$ be a wreath product of group $A = \langle X \rangle$ and $B = \langle Y \rangle$ such that $\pw(A, X) < \infty$ and $\pw(B, Y) < \infty$. Is it true that $\pw(G, X \cup Y) < \infty$?
\end{problem}

Riley and Sale \cite{RS} have proved that if $B$ is a finitely generated abelian group then the answer is positive. Fink \cite{fink} has also proved the finiteness of palindromic width of wreath product  for some more cases. However, the general case when $B$ is a finitely generated non-abelian group is still open.

For $g$, $h$ in $G$, the \emph{commutator} of $g$ and $h$ is defined as $[g, h]=g^{-1} h^{-1}gh$. If $\mathcal{C}$ is the set of commutators in some group $G$ then the commutator  subgroup $G'$ is generated by $\mathcal{C}$. The \emph{commutator length} $l_{\mathcal{C}}(g)$ of an element $g \in G'$ is the minimal number $k$ such that $g$ can be expressed as a product of $k$ commutators.  The \emph{commutator width} of $G$ is defined to be 
$$
{\rm cw}(G) = \underset{g \in G}{\sup} \ l_{\mathcal{\mathcal C}}(g). 
$$

It is well known \cite{R} that the commutator width of a free non-abelian group is infinite, but the commutator width of a finitely generated nilpotent group is finite (see \cite{AR,AR1}). An algorithm of the computation of the commutator length in free non-abelian groups can be found in \cite{B}.

Bardakov and  Gongopadhyay \cite[Problem 2]{BG} asked to understand a connection between commutator width and palindromic width. In \cite{BG3} they investigated this question and obtained further results in this direction. But still, the exact connection is subject of further inves\-ti\-ga\-tion.

Let us also note the following problem concerning palindromic width in solvable groups.

\begin{problem}[V. Bardakov~--- K. Gongopadhyay]
Is it true that the palindromic width of a finitely generated solvable group of step $r \geq 4$ is finite?
\end{problem}

The same problem is also open for the commutator width of solvable groups \cite[Question 4.34]{Kourovka}.

It is proved in \cite{BG3} that the palindromic width of a solvable Baumslag~--- Solitar group is equal 2. For the non-solvable groups  we formulate the following.

\begin{problem} [V. Bardakov~--- K. Gongopadhyay]
Let
$$
BS(m,n) = \langle a, t~||~t^{-1} a^m t = a^n \rangle,~~m, n \in \mathbb{Z} \setminus \{0\}
$$
be a non solvable Baumslag~--- Solitar group. Is it true that $\pw(BS(m,n), \{a, t \})$ is infinite?
\end{problem}

Let us formulate the general question on the palindromic width of HNN--extensions.

\begin{problem}[V. Bardakov~--- K. Gongopadhyay]
Let $G = \langle X \rangle$ be a  group and $A$ and $B$ are proper isomorphic subgroups of $G$ and $\varphi : A \rightarrow B$ be an isomorphism. Is it true that the HNN-extension
$$
G^* = \langle G, t~||~t^{-1} A t = B, \varphi \rangle
$$
of $G$ with associated subgroups $A$ and $B$ has infinite palindromic width with respect to the generating set $X \cup \{ t \}$?
\end{problem}

In \cite{BT} palindromic widths of free products were investigated. For the generalized free products we  formulate:

\begin{problem}[V. Bardakov~--- K. Gongopadhyay]
Let $G = A*_C B$ be a free product of $A$ and $B$ with amalgamated subgroup $C$ and $|A : C| \geq 3$, $|B : C| \geq 2$. Is it true that $\pw(G, \{A, B \})$ is infinite?
\end{problem}

Finally, we note the following problem raised by Riley and Sale.

\begin{problem}\cite{RS}
Is there a group $G$ with finite generating sets  $X$ and $Y$ such that ${\rm pw}(G, X)$ is finite but ${\rm pw}(G, Y)$ is infinite?
\end{problem}

\section{Brunnian words, Brunnian and other special braids} \label{sec2}

\textbf{\ref{sec2}.1 Brunnian words.}
Let $G$ be a group generated by a finite set $X$. An element $g\in G$ is called \emph{Brunnian} if there exists a word $w=w(X)$ on $X$ with $w=g$ such that, for each $x\in X$, $g=w$ becomes a trivial element in $G$ by replacing all entries $x$ in the word $w=w(X)$ to be $1$.

\smallskip

\noindent\textit{Example.} Let $X=\{x_1,\ldots,x_n\}$. Then $w=[[x_1,x_2],x_3,\ldots, x_n]$ is a Brunnian word in $G$. Any products of iterated commutators with their entries containing all elements from $X$ are Brunnian words.

\begin{problem}[J.~Wu]
 Given a group $G$ and a (finite) generating set $X$, find an algorithm  for detecting a Brunnian word that can NOT be given as a product of iterated commutators with their entries containing all elements from $X$.
\end{problem}

If $G$ is a free group with $X$ a basis, then all Brunnian words are given as products of iterated commutators with their entries containing all elements from $X$.

\smallskip

\noindent\textit{Most interesting Example.} Let $G=\langle x_0,x_1,\ldots,x_n \ | \ x_0x_1\cdots x_n=1\rangle$, where $X=\{x_0,x_1,\ldots,x_n\}$. For $n=2$, $[x_1,x_2]$ is Brunnian word in $G$. The solution of the question for this example may imply a combinatorial determination of homotopy groups of the $2$-sphere~\cite{Wu}.

\smallskip

\textbf{\ref{sec2}.2 Brunnian braids.}
The Brunnian braids over general surfaces have been studied in~\cite{BMVW}.

\begin{problem}{\rm (J.~Wu)}
\begin{enumerate}
\item {Find a basis for Brunnian braid group.}
\item {Determine Vassiliev Invariants for Brunnian braids. The relative Lie algebras of the Brunnian braid groups have been studied in~\cite{LVW}.}
\item {Classifying Brunnian links obtained from Brunnian braids.}
\item { Classifying links obtained from Cohen braids, where the Cohen braids were introduced in~\cite{BVW}.}
\end{enumerate}
\end{problem}

Note that Question 23 in \cite{Birman} asks to determine a basis for Brunnian braid group over $S^2$. A connection between Brunnian braid groups over $S^2$ and the homotopy groups $\pi_*(S^2)$ is given in~\cite{BCWW}.

\smallskip

\textbf{\ref{sec2}.3. Two-colored knot theory.} Let $Q$ be a subgroup of $P_n$. Then the free product $B_n\ast_{Q}B_n$ with amalgamation can be described as
$$
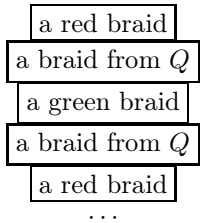
\begin{array}{c}
\fbox{\textrm{a red braid}}\\
\fbox{a braid from $Q$}\\
\fbox{a green braid}\\
\fbox{a braid from $Q$}\\
\fbox{a red braid}\\
\cdots\\
\end{array}
$$

\begin{problem}[J.~Wu]
Develop $2$-colored knot theory. What are the links obtained from $B_n\ast_QB_n$?
\end{problem}

There is a connection between the groups $P_n\ast_Q P_n$ and homotopy groups of spheres~\cite{MW}.

\smallskip

\textbf{\ref{sec2}.4. Finite type invariants of Gauss knots and Gauss braids.} Gibson and Ito \cite{GI}  defined  finite type invariants of Gauss knots.
It is not difficult to define finite type invariants of Gauss braids.

\begin{problem}[V. Bardakov]
 Is it true that  finite type invariants classify Gauss braids (Gauss knots)?
\end{problem}

\smallskip

\textbf{\ref{sec2}.5. Gr\"obner-Shirshov bases.}  About Gr\"obner-Shirshov bases see, for example, \cite{BC1, BC2}.

\begin{problem}[L. Bokut]
Find  Gr\"obner - Shirshov bases for generalizations of braid groups (Artin groups of types $B_n$, $C_n$, $D_n$, virtual braid groups and so on).
\end{problem}

\section{Braid groups, link groups, and 3-manifold groups} \label{sec3}

\textbf{\ref{sec3}.1. Intersecting subgroups of link groups.} Let $L_n$ be an $n$-component link. Let $R_i$ be the normal closure of the $i\,$th meridian.

\begin{problem}[J.~Wu]
Determine the group
$$
\frac{R_1\cap R_2\cap\cdots\cap R_n}{[R_1,R_2,\ldots, R_n]_S},
$$
where $[R_1,R_2,\ldots,R_n]_S=\prod\limits_{\sigma\in \Sigma_n}[[R_{\sigma(1)},R_{\sigma(2)}],\ldots, R_{\sigma(n)}]$.
\end{problem}

The cases $n=2,3$ have been discussed in~\cite{Wu2} (see also~\cite{LiWu}).  A link $L_n$ is called strongly nonsplittable if any nonempty proper sublink of $L_n$ is nonsplittable. In the case that $L_n$ is strongly nonsplittable, the group in the question is isomorphic to the homotopy group $\pi_n(S^3)$.

\smallskip

\textbf{\ref{sec3}.2. Extreme properties of hyperbolic 3-orbifold groups}. Fundamental groups of orientable hyperbolic 3-manifolds (in particular, knot complements) and groups of orientable hyperbolic 3-orbifolds have explicit presentations in the group ${\rm PSL}(2,\mathbb{C})$, isomorphic to the group of all orientation-preserving isometries of the hyperbolic 3-space. It is one of the key problems in the theory of hyperbolic 3-manifolds and 3-orbifolds to answer the question where a given subgroup of in ${\rm PSL}(2,\mathbb{C})$ is discrete. In 1977 J{\o}rgensen \cite{J1} proved that the question of a discreteness of arbitrary groups can be reduced to the question of discreteness of two-generated groups. 

In \cite{J2} J{\o}rgensen obtained a necessary condition for the discreteness of a two-generated group. The condition is presented by case (1) of theorem~\ref{theorem3}  and it looks as a nonstrict inequality connecting the trace of one of the generators and the trace of the commutator. Two more necessary discreteness conditions of similar form (see cases (2) and (3) of theorem~\ref{theorem3}) were later obtained by Tan~\cite{T1} and independently by Gehring and Martin~\cite{GM1}. Summarizing the mention results we have 

\begin{theorem} \label{theorem3}
Let $f, g \in {\rm PSL}(2,\mathbb{C})$  generate a discrete group. Then the properties hold: 
\begin{enumerate}
\item  If $\langle f,  g \rangle$ is nonelementary, then
$$
| \text{\rm tr}^{2} (f) - 4| + |\text{\rm tr} [f,g] - 2 | \geq 1. 
$$
\item If $\text{\rm tr} [f, g] \neq 1$, then 
$$
| \text{\rm tr}^{2} (f) - 2| + | \text{\rm tr} [f,g] - 1| \geq 1. 
$$
\item  If $\text{\rm tr}^{2} (f) \neq 1$, then
$$
| \text{\rm tr}^{2} (f) - 1| + | \text{\rm tr} [f,g]| \geq 1. 
$$ 
\end{enumerate}
\end{theorem} 

The nonelementary discrete groups, having such a pair of generators that the inequality in the case (1) becomes equality, are called \emph{J{\o}rgensen groups}. It was shown by Callahan~\cite{C1} that the figure-eight knot group is the only 3-manifold group which is a J{\o}rgensen group.  Analogously, the discrete groups, having such a pair of generators there the inequality in the case (2), respectively in the case (3), becomes equality, are called \emph{Tan groups}, respectively, \emph{Gehring~--- Martin~--- Tan groups}. It was shown by Masley and Vesnin \cite{VM1} that the group of the figure-eight orbifold with singularity of order four is a Gehring~--- Martin~--- Tan group. 

\begin{problem}[A.~Vesnin]
Find all hyperbolic 3-orbifold groups which are Gehring~--- Martin~--- Tan groups or Tan groups. 
\end{problem}

\smallskip 

\textbf{\ref{sec3}.3. Cabling operations on some subgroups of ${\rm PSL}(2,\mathbb{C})$.}  Let $K$ be a framed knot with its group as a subgroup of ${\rm PSL}(2,\mathbb{C})$. Let $K_n$ be the $n$-component link obtained by naive cabling of $K$ along the frame.

\begin{problem}[J.~Wu]
Describe a cabling construction for the link group $\pi_{1}(S^{3} \setminus  K_n)$ as a subgroup of ${\rm PSL}(2,\mathbb{C})$.
\end{problem}

The groups of $K_n$ with $n\geq1$ admit a canonical simplicial group structure with its geometric realization homotopic to the loop space of $S^3$ if $K$ is a nontrivial framed knot (see~\cite{LLW}).

\smallskip

\textbf{\ref{sec3}.4. Linearity of 3-manifold groups.} Let $M^3$ be a 3-manifold. It is known that the fundamental group $\pi_1(M^3)$ is residually finite. Also, if $M^3$ is a standard manifold (Nil-manifold, Sol-manifold, hyperbolic and so on) then the fundamental group $\pi_1(M^3)$ is linear.

\begin{problem}[V. Bardakov]
Is it true that the group $\pi_1(M^3)$ is linear in general case? In particular, is it true that this group is linear for the case then $M^3$ is the compliment of some link $L$ in 3-sphere, i.e. $M^3 = S^3 \setminus L$?
\end{problem}

Every knot $K$ in $S^3$ is a torus knot, hyperbolic knot or satellite knot. Groups of torus knots and hyperbolic knots are linear. Hence, it is need to study only satellite knots. For questions related to the above problem and much more on 3-manifold groups we refer to the recent survey \cite{aus}. 

\smallskip

\textbf{\ref{sec3}.5. Inner and quasi-inner automorphisms.}  M. Neshchadim proved that for the group $G = \langle a, b ~||~ a^2 = b^2 \rangle$ the group $\mathrm{Aut}_c(G)$ of class-preserving automorphisms  is not equal to the group $\mathrm{Inn}(G)$ of inner automorphisms. This group $G$ is a free product of two cyclic groups with amalgamation $G = \mathbb{Z} *_{2\mathbb{Z} = 2\mathbb{Z}} \mathbb{Z}$ and is a group of torus link. On the over hand the group of trefoil knot is isomorphic to the braid group $B_3 = \mathbb{Z} *_{3\mathbb{Z} = 2\mathbb{Z}} \mathbb{Z}$ and for this group $\mathrm{Aut}_c(B_3) = \mathrm{Inn}(B_3)$. Also we know that if $K$ is alternating knot, then its group $\pi_{1} (S^{3} \setminus K)$ has a decomposition  $\pi_{1} (S^{3} \setminus K) = F_n *_{F_{2n-1}} F_n$ for some natural $n$.

\begin{problem}[V. Bardakov]
For which knot (link) $K$ its group $\pi_{1} (S^{3} \setminus K)$ has the property:
$$
\mathrm{Aut}_c(\pi_{1} (S^{3} \setminus K)) = \mathrm{Inn}(\pi_{1} (S^{3} \setminus K))?
$$
\end{problem}

\smallskip

\textbf{\ref{sec3}.6. Free nilpotent groups.} Consider a  reduced word $w = w(a,b)$ in the free group $F_2 = \langle a, b \rangle$ such that the product $x * y = w(x, y)$ defines a group operation in the free nilpotent group $F_n / \gamma_k (F_n)$ for all $n \geq 2$ and $k \geq 1$.

\begin{problem}[M. Neshchadim]
 Find all words $w$ having above property. Is it true that $w = a b$ or $w = b a$?
\end{problem}

\section{Symmetric cohomology of groups} \label{sec4}

\textbf{\ref{sec4}.1. Cohomology of groups.} Cohomology of groups is a contravariant functor turning groups and modules over groups into graded abelian groups. It came into being with the fundamental work of Eilenberg and MacLane \cite{Eilenberg1,Eilenberg2}. The theory was further developed by Hopf, Eckmann, Segal, Serre, and many other mathematicians. It has been studied from different perspectives with applications in various areas of mathematics, and provides a beautiful link between algebra and topology.

We recall the construction of the cochain complex defining the cohomology of groups. Let $G$ be a group and $A$ a $G$-module. As usual $A$ is written additively and $G$ is written multiplicatively, unless otherwise stated. For each $n \geq 0$, let $C^n(G, A)$ be the group of all maps $\sigma: G^n \to A$. The coboundary map $$\partial^n:C^n(G, A) \to C^{n+1}(G, A)$$ is defined as
\begin{equation*}\label{eqn1}
\begin{split}
\partial^n(\sigma)(g_1,\dots,g_{n+1})& = g_1\sigma(g_2,\dots,g_{n+1}) \\
& +  \sum_{i=1}^n (-1)^{i}\sigma(g_1,\dots,g_ig_{i+1},\dots,g_{n+1})  +(-1)^{n+1} \sigma(g_1,\dots,g_n),
\end{split}
\end{equation*}
for $\sigma \in C^n(G, A)$ and $(g_1,\dots,g_{n+1}) \in G^{n+1}$. Since $\partial^{n+1}\partial^n=0$, $\{C^*(G, A), \partial^*\}$ is a cochain complex. Let $Z^n(G, A)= {\rm Ker} (\partial^n)$ and $B^n(G, A)= {\rm Image} (\partial^{n-1})$. Then the $n$th cohomology of $G$ with coefficients in $A$ is defined as $$H^n(G,A)= Z^n(G, A)/B^n(G, A).$$

Cohomology of groups has concrete group theoretic interpretations in low dimen\-sions. More precisely, $H^0(G,A) = A^G$ and $Z^1(G,A) =$ Group of derivations from $G$ to $A$. Let $\mathcal{E}(G,A)$ be the set of equivalence classes of extensions of $G$ by $A$ giving rise to the given action of $G$ on $A$. Then there is a one-one correspondence between $H^2(G,A)$ and $\mathcal{E}(G,A)$. There are also group theoretic interpretations of the functors $H^n$ for $n \geq 3$.

Let $\Phi:H^2(G,A) \to \mathcal{E}(G,A)$ be the one-one correspondence mentioned above. Under $\Phi$, the trivial element of $H^2(G,A)$ corresponds to the equivalence class of an extension $$0 \to A \to E \to G \to 1$$ admitting a section $s:G \to E$ which is a group homomorphism.

An extension $\mathcal{E}:0 \to A \to E \to G \to 1$ of $G$ by $A$ is called a symmetric extension if there exists a section $s:G \to E$ such that $s(g^{-1})=s(g)^{-1}$ for all $g \in G$. Such a section is called a symmetric section. Let $\mathcal{S}(G,A)=\big\{[\mathcal{E}] \in \mathcal{E}(G,A)~|~\mathcal{E}~\textrm{is a symmetric extension}\big\}$. Then the following question seems natural.

\begin{question} \label{quest1}
What elements of $H^2(G,A)$ corresponds to $\mathcal{S}(G,A)$ under $\Phi$?
\end{question}

Staic \cite{Staic1, Staic2} answered the above question for (abstract) groups. Motivated by some questions regarding construction of invariants of 3-manifolds, Staic introdu\-ced a new cohomology theory of groups, called the symmetric cohomology, which classifies symmetric extensions in dimension two.

\smallskip

\textbf{\ref{sec4}.2. Symmetric cohomology of groups.} In this section, we present Staic's construction of symmetric cohomology of groups. For topological aspects of this construction, we refer the reader to \cite{Staic1, Staic2}. For each $n \geq 0$, let $\Sigma_{n+1}$ be the symmetric group on $n+1$ symbols. In \cite{Staic1}, Staic defined an action of $\Sigma_{n+1}$ on  $C^n(G, A)$. Since the transpositions of adjacent elements form a generating set for $\Sigma_{n+1}$, it is enough to define the action of these transpositions $\tau_i = (i, i+1)$ for $1 \leq i \leq n$. For $\sigma \in C^n(G,A)$ and $(g_1,...,g_n) \in G^n$, define
\begin{equation*}
\begin{split}
(\tau_1 \sigma)(g_1,g_2,g_3, \dots,g_n) &= -g_1 \sigma\big(g_1^{-1},g_1g_2,g_3,\dots,g_n\big),\\
(\tau_i \sigma)(g_1,g_2,g_3,\dots,g_n) &= -\sigma\big(g_1,\dots,g_{i-2},g_{i-1}g_i, g_i^{-1}, g_ig_{i+1}, g_{i+2},\dots,g_n\big) \\  
& \qquad \qquad \qquad \qquad \qquad \qquad \qquad \qquad \qquad \qquad  \textrm{for}~ 1 < i < n,\\
(\tau_n \sigma)(g_1,g_2,g_3,\dots,g_n) &= -\sigma\big(g_1,g_2,g_3,\dots,g_{n-1}g_n, g_n^{-1}\big).
\end{split}
\end{equation*}

It is easy to see that
\begin{equation*}
\begin{split}
\tau_i\big(\tau_i (\sigma)\big) &=\sigma,\\
\tau_i \big(\tau_j(\sigma) \big) &= \tau_j \big(\tau_i(\sigma) \big)~ \text{for} ~j \neq i\pm 1,\\
\tau_i\big(\tau_{i+1}\big(\tau_i(\sigma)\big)\big) &= \tau_{i+1}\big(\tau_i\big(\tau_{i+1}(\sigma)\big)\big).
\end{split}
\end{equation*}

Thus there is an action of $\Sigma_{n+1}$ on $C^n(G, A)$. Define $d^j: C^n(G,A) \to C^{n+1}(G,A)$ by
\begin{equation*}
\begin{split}
d^0(\sigma)(g_1,\dots,g_{n+1}) & =g_1\sigma(g_2,\dots,g_{n+1}),\\
d^j(\sigma)(g_1,\dots,g_{n+1}) & =\sigma(g_1,\dots, g_jg_{j+1}, \dots,g_{n+1})~ \text{for} ~1 \leq j \leq n,\\
d^{n+1}(\sigma)(g_1,\dots,g_{n+1}) & =\sigma(g_1,\dots ,g_n).
\end{split}
\end{equation*}

Then $\partial^n(\sigma)= \sum_{j=0}^{n+1} (-1)^j d^j(\sigma)$. It turns out that
\begin{equation*}\label{eqn2}
\begin{split}
\tau_id^j &=d^j \tau_i~\text{if}~i<j,\\
\tau_id^j &=d^j \tau_{i-1}~\text{if}~j+2 \leq i,\\
\tau_id^{i-1} & =-d^i,\\
\tau_id^i & =-d^{i-1}.
\end{split}
\end{equation*}

Let $CS^n(G, A)= C^n(G, A)^{\Sigma_{n+1}}$ be the group of invariant $n$-cochains. If $\sigma \in CS^n(G, A)$, then it follows from the above identities that $\partial^n(\sigma) \in CS^{n+1}(G, A)$.
Thus the action is compatible with the coboundary operators and we obtain a cochain complex $\{CS^*(G, A), \partial^*\}$. Its cohomology, denoted by $HS^*(G,A)$, is called the symmetric cohomology of $G$ with coefficients in $A$. In \cite {Staic1, Staic2}, Staic gave examples of groups for which the symmetric cohomology is different from the ordinary cohomology.

As in ordinary cohomology, $HS^0(G,A)=A^G$ and $ZS^1(G,A)=$ Group of symmetric derivations from $G$ to $A$. The inclusion $CS^*(G, A) \hookrightarrow C^*(G, A)$ induces a homomorphism $h^*:HS^*(G, A) \to H^*(G, A)$. In \cite{Staic2}, Staic proved that the map $$h^*:HS^2(G, A) \to H^2(G, A)$$ is injective and the composite map $$\Phi \circ h^*: HS^2(G, A) \to \mathcal{S}(G,A)$$ is bijective. Thus the symmetric cohomology in dimension two classifies symmetric extensions. This answers Question~\ref{quest1}.

When a group under consideration is equipped with a topology or any other structure, it is natural to look for a cohomology theory which also takes the topology or the other structure into account. This lead to various cohomology theories of topological groups and Lie groups. Topology was first inserted in the formal definition of cohomology of topological groups in the works of S.-T. Hu \cite{Hu}, W.T. van Est \cite{vanEst} and A. Heller \cite{Heller}. In a recent paper \cite{Singh}, Singh studied continuous and smooth versions of Staic's symmetric cohomology. He defined a symmetric continuous cohomology of topological groups and gave a characterization of topolo\-gi\-cal group extensions that correspond to elements of the second symmetric con\-ti\-nuous cohomology. He also defined symmetric smooth cohomology of Lie groups and proved similar results in the smooth setting. His results also answered continuous and smooth analogous of Question~\ref{quest1}.

\smallskip

\textbf{\ref{sec4}.3. Some Problems.} Symmetric cohomology of groups is a fairly recent construc\-tion, and not much is known about its theoretical and computational aspects. In view of this, the following questions seems natural. We hope that answers to these questions will help in better understanding and possible applications of symmetric cohomology.

\begin{problem}[M. Singh]
Does there exists a Lyndon~--- Hochschild~---- Serre type spectral sequence for symmetric cohomology of groups?
 \end{problem}

\begin{problem}[M. Singh]
Cohomology of a discrete group can also be defined as the cohomology of its classifying space. Is there a topological way of obtaining symmetric cohomology of discrete groups?
 \end{problem}

\begin{problem}[M. Singh]
Is it possible to define a symmetric cohomology of Lie algebras? How does this relate to the symmetric cohomology of Lie groups?  It seems possible to do so, and Staic suspect that it is equal to the usual cohomology \cite{Staic3}.
 \end{problem}

\begin{problem}[M. Singh]
Fiedorowicz and Loday \cite{Fiedorowicz} introduced crossed simplicial groups and a cohomology theory of these objects, which is similar to symmetric cohomology. It would be interesting to explore some connection between the two~\cite{Staic1}.
 \end{problem}

\begin{problem}[M. Singh]
Let $G$ be a group and $\mathbb{C}^{\times}$ a trivial $G$-module. Then the Schur multiplier of $G$ is defined as $\mathcal{M}(G):=H^2(G, \mathbb{C}^{\times})$. It turns out that the Schur multiplier $\mathcal{M}(G)$ of a finite group $G$ is a finite abelian group. Finding bounds on the order of $\mathcal{M}(G)$ is an active area of research and has wide range of applications, particularly in automorphisms and representations of finite groups. We define the symmetric Schur multiplier of $G$ as $$\mathcal{M}S(G):= HS^2(G, \mathbb{C}^{\times}).$$ Clearly, the symmetric Schur multiplier $\mathcal{M}S(G)$ of a finite group $G$ is a finite group. It would be interesting to find bounds on the order of $\mathcal{M}S(G)$ and $\mathcal{M}(G)/\mathcal{M}S(G)$.
 \end{problem}

\begin{problem}[M. Singh]
Let $G$ be a finite group. Then the Bogomolov multiplier of $G$ is defined as $$B_0(G) = {\rm ker} \big(res^G_A: H^2(G, \mathbb{C}^{\times}) \to \bigoplus_{A \subset G} H^2(A, \mathbb{C}^{\times}) \big)$$ where $A$ runs over all abelian subgroups of $G$ and $res^G_A$ is the usual restriction homomorphism. The group $B_0(G)$ is a subgroup of the Schur multiplier $\mathcal{M}(G)$, and appears in classical Noether's problem and birational geometry of quotient spaces of $G$. See \cite{Kang} for a recent survey article. It would be interesting to find relation between $\mathcal{M}S(G)$ and $B_0(G)$.
 \end{problem}


\end{document}